\documentclass[a4paper,10pt]{article}
\usepackage{amsfonts}
\usepackage{amssymb}
\usepackage[utf8]{inputenc}

\usepackage[english ]{babel}

\usepackage{bbding}

\usepackage[all]{xy}

\usepackage{tikz-cd} 
\usepackage[all]{xy}
\usepackage{amsfonts}
\usepackage{amssymb}
\usepackage{graphicx}
\usepackage{mathtools}
\usepackage{skak} 
\usepackage{foekfont}

\usepackage{calligra}
\usepackage[T1]{fontenc}

\usepackage{amsmath,mathdots}

\usepackage{stmaryrd}

\usepackage[utf8]{inputenc}

\usepackage{mathtools}
\usepackage{mathdots}
\usepackage{tikz}

\definecolor{amber(sae/ece)}{rgb}{1.0, 0.49, 0.0}
\newfont{\rsfsten}{rsfs10 scaled 1200}
\usepackage{MnSymbol}
\usepackage{tikz}
\usepackage{amscd}
\usepackage{MnSymbol}
\usepackage{wasysym}

\usepackage{mathrsfs}

\usepackage{pdfcolmk}

\usepackage{graphicx}

\newcommand*{\rom}[1]{\expandafter\@slowromancap\romannumeral #1@}

\usepackage{wrapfig}

\newcommand{\tightunderset}[2]{%
  \mathop{#2}\limits_{\vbox to .3ex{\kern-0.95ex\hbox{$#1$}\vss}}}
\newcommand{\tightoverset}[2]
{%
  \mathop{#2}\limits_{\vbox to .3ex{\kern-0.95ex\hbox{$#1$}\vss}}}

\newcommand{\oset}[2]{%
  {\mathop{#2}\limits^{\vbox to -.5\ex@{\kern-\tw@\ex@
   \hbox{\scriptsize #1}\vss}}}}
\usepackage[utf8]{inputenc}

\usepackage {ifsym}

\usepackage {url}

\makeatother

\title{  A Dozen  Problems, Questions and Conjectures  about Positive Scalar Curvature}
\author{Misha Gromov}

\begin{document}

\maketitle \tableofcontents

\begin{abstract} 

Unlike manifolds with positive   sectional and with positive  Ricci curvatures which aggregate to  modest (roughly) convex   islands in the vastness of all  Riemannian spaces, the domain  $\{{\cal SC}>0\}$ of manifolds with  {\it positive scalar curvatures}  protrudes in all direction as a gigantic octopus or an  enormous  multi-branched  tree. 
  Yet, there are  certain   rules to the shape of $\{{\cal SC}>0\}$    which limit the spread of this domain but most of these rules remain a guesswork.

  In the  present  paper  we collect a few "guesses"  extracted from  a  longer  article, which is still in preparation: {\sl 100 Questions, Problems and Conjectures around  Scalar Curvature.}
  
   Some of these "guesses" are presented as {\sf questions} and some  as {\sf conjectures}. Our  formulation of  these conjectures  is not supposed to be either   most general or    most plausible, but rather   maximally thought provoking.

\end{abstract}

\section {Definition of   Scalar Curvature.} 
  
 The scalar curvature of a $C^2$-smooth  Riemannian manifold $X=(X,g)$, denoted 
 $Sc=Sc(X)=Sc(X,g)=Sc(g)$ is a continous  function on $X$, written  as $Sc(X)(x)$ and  $Sc(g)(x)$,
 $x\in X$, 
 which is uniquely characterised by the following four  properties.

 $\bullet_1$  {\it  Additivity under Cartesian-Riemannian Products .} 
 $$Sc(X_1\times X_2, g_1\oplus g_2)=Sc(X_1,g_1)+ Sc(X_2, g_2),$$ 
 where this equality is understood point-wise,
 $$Sc(X_1\times X_2)(x_1,x_2) =Sc(X_1)(x_1)+ Sc(X_2)(x_2).$$

   $\bullet_2$  {\it  Scale covariance. }
   $$Sc( X, \lambda^2\cdot g)=\lambda^2\cdot Sc(X)\mbox { for all real } \lambda>0.$$
   
Thus, for instance, since   $(\mathbb R^n, g_0) $  is isometric to  $(\mathbb R^n, \lambda^2\cdot g_0)$ 
for the Euclidean metric $g_0$, 
$$Sc(\mathbb R^n)=0\mbox  {  for all } n=1,2,3,... .$$
 
 \vspace {1mm}

     $\bullet_3$  {\it Volume Comparison.} If the scalar curvatures of $n$-dimensional manifolds $X$ and $X'$ at some points $x\in X$  and  $x'\in X'$ are related by the strict inequality
  $$Sc(X)(x)< Sc(X')(x'),$$
then the Riemannian volumes of  small balls around these points   satisfy
  $$vol(B_{x}(X, \varepsilon)) > vol(B_{x'}(X', \varepsilon))$$
  for all sufficiently small $\varepsilon>0$.

 This volume inequality, in agreement with  $\bullet_1$, is  {\it additive under Riemannian products}:
 if 
 $$\mbox {$vol(B_{x_i}(X, \varepsilon)) > vol(B_{x'_i}(X_i', \varepsilon))$, for $\varepsilon \leq \varepsilon_0$,}$$ 
and {\it for  all  points}   $x_i\in  X_i$ and $x'_l\in X'_i$,  $i=1,2$,
then 
 $$ vol_n(B_{(x_1, x_2)}(X_1\times X_2, \varepsilon_0)) > vol_n(B_{(x'_1,x_2')}(X'_1\times X'_2, \varepsilon_0)$$
   {\it for  all}   $(x_1,x_2) \in  X_i\times X_2 $ and  $(x'_1,x'_2) \in  X'_1\times X'_2 $.
 
 This follows from the Pythagorean  formula 
 $$dist_{X_1\times X_2} =\sqrt{dist_{X_1}^2+dist_{X_2}^2}.$$ 
 and  the Fubini theorem applied to the  "fibrations" of balls over balls:
 $$B_{(x_1, x_2)}(X_1\times X_2, \varepsilon_0)) \to B_{x_1}(X_1,\varepsilon_0)\mbox  { and } 
B_{(x'_1, x'_2)}(X'_1\times X'_2, \varepsilon_0))  \to B_{x_1}(X'_1,\varepsilon_0),$$
 where the fibers are balls of radii $\varepsilon \in [0,\varepsilon_0]$ in $X_2$ and $X'_2$.
 \vspace {1mm}
 
  $\bullet_4$  {\it Normalisation/Convention  for $2$-spheres.} The unit sphere $S^2=S^2(1)$ has constant scalar curvature  $2$ (twice the sectional curvature).  \vspace {1mm}
 
It is an elementary exercise to prove the following. \vspace {1mm}

{\Large $\mathbf \star$}$_1$  {\sl The function $Sc(X,g)(x)$ which satisfies $\bullet_1$-$\bullet_4$   exists  and unique};\vspace {1mm}
  
{\Large $\mathbf \star$}$_2$  {\sl  The unit spheres and  the hyperbolic spaces with $sect.curv=-1$ satisfy
  $$Sc(S^n(1))=n(n-1)\mbox  {  and } Sc(H^n(-1))=-n(n-1).$$}
  Thus, 
  $$Sc(S^n(1)\times H^n(-1))=0=Sc(\mathbb R^n),$$
 which implies  that the volumes of the  small balls in $S^n(1)\times H^n(-1)$ are "very close" to the volumes of the Euclidean $2n$-balls.
 \vspace {1mm}

{\Large $\mathbf \star$}$_3$   {\sl The scalar curvature of a Riemannin manifold $X$ is equal to the sum of the values  of the   sectional curvatures  at the bivectors  of  an orthonormal frame\footnote {Remarkably, this sum is independent of the  frame by the Pythagorean theorem.}  in $X$, 
 $$Sc(X)(x)=\sum_{i, j} c_{ij},\mbox {  }  i,j=1,...,n.$$}\vspace {1mm}

For example, all {\it compact Riemannin symmetric spaces} $X$, except for the $n$-torus $\mathbb T^n$, have $Sc(X)>0$,  
while $\mathbb T^n$, being covered by $\mathbb R^n$, has  $Sc (\mathbb T^n)=0$.

\vspace {1mm}
 
It may be tempting to take  the above $\bullet_1-\bullet_4$ for a definition of scalar curvature for singular metric spaces $X$.  In fact, it may work  for  $X$ with moderate singularities, e.g. for {\it Alexandrov's spaces with sectional  curvatures bounded from below}  (see [1]\footnote {[1]  S.Alexander, V. Kapovitch, A.Petrunin,
{\sl Alexandrov 
geometry,}

{\url http://www.math.psu.edu/petrunin/}}), where the properties of the  so defined scalar curvature  must be comparable  to what is observed  in the smooth case (see section 7).

Yet, volumes of balls to not touch the heart of the scalar curvature;  we suggests an   alternative in section 7.

\section {Soft and Hard Facets of Scalar Curvature.}

  We are not so much concerned with the scalar curvature $Sc(X)$ per se, but rather  with the   
  effect of {\it lower scalar curvature bounds} on  the geometry and the topology of $X$, where, for instance,  the inequality    "$Sc(X)>0$" 
   can be   \textbf {defined}  by saying that   \vspace {1mm}
 
 {\it all sufficiently small balls $B_x(\varepsilon)\subset  X$, $\varepsilon \leq \varepsilon_0(x)>0$, have the volumes smaller than the volumes of the equividimensional  Euclidean $\varepsilon$-balls.} 

  \vspace {1mm}
Then   "$Sc(X)\geq 0$" is defined as 

  \vspace {1mm}

\hspace {30mm} {\it $Sc(X)>-\varepsilon  $"  for all $\varepsilon >0$.}

  \vspace {1mm}

Similarly   \vspace {1mm}

{\it "$Sc(X)\geq \sigma$", $\sigma >0$, is equivalent  the volumes of $B_x(\varepsilon)$ in $X$   being smaller than the volumes of the  $\varepsilon$-balls in the Euclidean  spheres $S^n(R)$ of radii $R>\sqrt{(n(n-1)/\sigma)}$}, \vspace {1mm}

\hspace {-6mm} and  $Sc(X)\geq -\sigma$ is expressed by   \vspace {1mm}

{\it the  bound on the   volumes of $B_x(\varepsilon)$    by those of the $\varepsilon $-balls in the hyperbolic spaces with constant  the sectional curvatures $<-\sigma/n(n-1)$. } \vspace {1mm}

Alternatively, "$Sc(X)\geq -\sigma$" can be defined with  {\it no reference} to hyperbolic spaces by {\it the reduction to the case $\sigma=0$} and appealing to the relation
$$Sc(X\times S^m(R))\geq 0
\mbox{  for    $R=\sqrt{(m(m-1)/\sigma}$},$$ where one may use any $m\geq 2$ one likes.
\vspace {1mm}

 Although  the key role of the  scalar curvature in general relativity  was  established by  Hilbert's 
 variational derivation of the Einstein equation more than  a century ago (see [2]\footnote {[2] David Hilbert, The Foundations of Physics", 1915.})  the significance of $Sc(X)$ in the global geometry and in topology remained obscure until 1963, when Andr\'e    Lichnerowicz (see [3]\footnote {[3]  A. Lichnerowicz,
Spineurs harmoniques.
C. R. Acad. Sci. Paris, S\'erie A, 257 (1963), 7-9.}) showed    that the inequality  $Sc(X)>0$  imposes non-trivial constraints on the topology of $X$.

For instance,    Lichnerowicz' theorem implies that  \vspace {1mm}

{\it if $m$ is even, then  smooth    complex  projective hypersurfaces $X\subset \mathbb CP^{m+1}$}   (these have  real dimension $ dim(X)=2m$)  {\sl  of degrees $\geq m+2$,  e.g. $X\subset \mathbb CP^3$  given by the equation 
  $$x^4_1+x^4_2+x^4_3+x^4_4=0,$$
 admit  {\large no}  metrics with $Sc>0$.}\vspace {1mm}

This follows from the   Atiyah-Singer formula for the (Atiyah-Singer)-Dirac operator  $D$ confronted with (what is now called) the {\it Schroedinger-Lichnerowicz-(Weitzenboeck-Bochner) identity.}

In fact, the index formula  implies that the  index of $D$ on these  manifolds {\it does not vanish},\footnote{This formula says in the present case that $Ind(D)=\hat A(X)$ where $\hat A(X)$ is a  particular Pontryagin number of $X$.} and, consequently,  {\it there are non-zero  harmonic spinors on these $X$} (i.e. solutions $s$ of $D(s)=0)$, while    the    {\it Schroedinger-Lichnerowicz-(Weitzenboeck-Bochner) identity}  
$$D^2=\nabla\nabla^\ast+\frac{1}{4} Sc,$$
 shows that closed manifolds  with $Sc>0$ admit {\it no harmonic spinors}. \vspace {1mm}

Eleven years  later, Nigel Hitchin(see [4]\footnote {[4] N. Hitchin,
Harmonic Spinors
, Adv. in Math.
14
(1974), 1-55.}) used a more sophisticated 1971  version of  the Atiyah-Singer
 index theorem which yields harmonic spinors on some {\it exotic spheres}  $\Sigma^n$ (which are homeomorphic but  not diffeomorphic to the ordinary  spheres $S^n$) of dimensions $n=8k+1$ and $n=8k+2$ and which, together  with the 
Schroedinger-Lichnerowicz' identity, implies that  

\hspace {25mm} {\sl there is no  metrics with $Sc>0$ on these  $\Sigma^n$.}\vspace {1mm}

Then Stefan  Stolz,  elaborating on  the earlier work by several authors, 
 showed that there are   

{\sl no further    obstructions to the existence  of metrics with $Sc>0$ on {\large \sf simply connected} manifolds  of dimension $\geq 5$ besides  those delivered by the index theorem} [5].\footnote {[5]  S. Stolz. {\sl Simply   connected  manifolds   of positive  scalar  curvature},
  Ann.  of  Math.  (2) 136
  (1992),  511-540.}

For instance 

{\sl all simply connected manifolds of dimensions $n= 3, 5, 6, 7 \mod 8$ admit metrics with positive scalar curvatures.} \vspace {1mm}

The proof of this theorem,  which  relies  on {\it surgery of manifolds with} $Sc>0$   and on the  {\it cobordism theory},  suggests that  manifolds with positive scalar curvature are  {\sf almost}  as  {\it soft}  as  smooth manifolds with no geometric constraints imposed  on them. 
But  the grand picture of   scalar curvature in all its beauty unravels  when one   looks beyond this    {\sf "almost"}. 

(The opposite inequality $Sc(X)<0$   is truly and fully    soft  and,  unlike  $Sc>0$, has no influence on the topology and global geometry of $X$ what-so-ever (see [6]\footnote {[6]   J. Lohkamp,
{\sl  Metrics of negative Ricci curvature}, Annals of Mathematics, 140 (1994), 655-683.}).
\vspace {1mm}
 
A manifestly     {\it rigid}  property of $Sc>0$  can be already  seen  in the following corollary to Schoen-Yau solution  
of the Riemannin positive mass conjecture in relativity (see[7]\footnote {[7]  R. Schoen and S.-T. Yau, On the proof of the positive mass conjecture in general relativity, Commun. Math. Phys. 65,  (1979). 45-76.}).

\vspace {1mm}

\textbf {Solution of the Geroch Conjecture.}\footnote {Attribution of this simplified podescribessitive mass conjecture  to Robert Geroch  is made in the above cited paper by Schoen and Yau.

In fact, the full  Riemannin positive mass conjecture which describes possible asymptotic   behaviours of metrics 
with $Sc>0$ on $\mathbb R^3$  (and on $\mathbb R^n$ for this matter) which are close  (rather than equal) to the Euclidean metric at infinity follows from this Geroch conjecture according to   {\sl J. Lohkamp, Scalar curvature and hammocks, Math. Ann. 313 (1999), 385-407}.}
  {\sl The Euclidean metric $g_0$  on   $\mathbb R^3$  (which has  $Sc(g_0)=0$) admits no non-trivial  compactly supported  perturbations $g$  with $Sc(g)\geq 0$. 

Namely, if  a smooth Riemannin  metric $g$ on  the Euclidean space $\mathbb R^3$ has  $Sc(g)\geq 0$ and if $g$ is equal to $g_0$ outside a compact subset in $\mathbb R^3$, then $Sc(g)=0$; moreover,  $g$ is  Riemannin flat, that is    $(\mathbb R^3, g)$ is isometric to $(\mathbb R^3,g_0)$}.    \vspace {1mm}

\vspace {1mm}

This result has been refined and generalised in a variety of directions (see below and also [13] and [21] at the end of the next section and references therein)
but 
{\large \sf  the rigidity of $Sc>0$}  we are after, albeit related to the above, is of different nature.  In fact what we look for is  \vspace {1mm}

{\it  a  structurally organised set of} ({\sf desirably sharp}) {\it geometric inequalities satisfied  by manifolds with $Sc >0$, more generally, with $Sc\geq \sigma$.  

Also, we
search for a   general category (or categories)   of  spaces, or other kind of  objects, which would satisfy (certain classes of)  such inequalities.}

\vspace {1mm}
 
\hspace {25mm}  {\large \sf \textbf {Additional Remarks and  References.}} \vspace {1mm}

Geroch conjecture has been validated in all dimensions:   \vspace {1mm}

{\sl The Euclidean metrics on $\mathbb R^n$  for all $n$ admit no non-trivial   compactly supported perturbations with
$Sc\geq0$.} \vspace {1mm}

 This (trivially) follows, for instance,   from {\it non-existence of metrics with $Sc>0$  on the $n$-tori}  where the latter  can be most easily proved by applying the index  theorem to  suitably "twisted" Dirac operators.  
\vspace {1mm}

Witten suggested  a different way of using the Dirac operator in the context of the positive mass problem,   where the index theorem is replaced by a direct proof of  {\sl harmonic stability of parallel spinors} on $\mathbb R^n$  under certain perturbations of the Euclidean metric.

By a similar method,  Min-Oo (see [8]\footnote{[8]   M. Min-Oo, Scalar curvature rigidity of asymptotically hyperbolic spin manifolds,. Math. Ann. 285, 527?539 (1989)}) proved that \vspace {1mm}

{\it the hyperbolic metric $g_0$ on the real hyperbolic space  $H_\mathbb R ^n$   admits non non-trivial compactly supported  perturbations $g$   with $Sc(g)\geq -n(n-1)=Sc(g_0)$.}\vspace {1mm}

Apparently, it is unknown \vspace {1mm}

{\large \sf if other symmetric spaces of non-compact types admit  compactly supported perturbations of their Riemannian metrics which would increase scalar curvature.}

\section {Bounds on the Uryson Width,  Slicing Area  and Filling Radius.}

 {\large \sf  \textbf {A.}} {\large{\sf Conjecture.  }} {\sl Let $X$ be an $n$-dimensional  Riemannin manifold with scalar curvature bounded from below by
 $$Sc (X)\geq n(n-1)=Sc(S^n).$$
Then the $(n-1)$-dimensional {\sf Uryson width} of $X$ is bounded by a universal constant. 

This means that 
 there exists a continuous map  from $X$ to an $(n-1)$-dimensional polyhedral space  $P$,
 $$f:X\to P=P^{n-1}, $$
 such that the pullbacks of all points  have controllably  bounded diameters, namely,
 $$diam_X(f^{-1}(p))\leq const\mbox  { for all } p\in P.$$
 for some universal constant $const>0$ possibly (and undesirably)  depending on $n$ .}\vspace {1mm}
 
This conjecture  says, in effect, that 
  that $n$-dimensional  manifolds $X$ with $Sc(X)\geq \sigma>0$   "topologically spread" in  at most $n-1$ directions. 
 
  In fact, one expects that these $X$
 spread only in $n-2$ direction which can be formulated as follows. \vspace {1mm}

 {\large \sf  \textbf {A$_+$}.} {\large{\sf Conjecture.}}{\sl The above $X$ admits a  a continuous map $f$ to an $(n-2)$-dimensional polyhedral space  $P$,
  such that 
 $diam_X(f^{-1}(p))\leq const_+$  for all $p\in P$.}
 \vspace {1mm}

 But the most attractive (and  least tenable) is the  conjecture  {\large \sf   {A$_{++}$}} below which claims  that {\it closed} manifolds  with $Sc\geq \sigma>0$ can be {\it sliced by surfaces with small areas} according the following definition. \vspace {1mm}

 {\large \sf Slicings and  Waists.} 
  An  {\it $m$-sliced $n$-cycle},  $m\leq n$, is an   $n$-dimensional  psedomanifold $P=P^n$  partitioned into {\it $m$-slices $P_q\subset P$}, which are  the pullbacks of the points of  a simplicial map $\varphi:P\to Q$ where $Q$ is an $(n-m)$-dimensional   pseudomanifold  and where all  pullbacks $P_q=\varphi^{-1}(q)\subset P$ have $dim(P_q)\leq m$, $q\in Q$. 
  
   (Sometimes one insists that  $\varphi$  must be   {\it proper}, hence, with compact pullbacks $\varphi^{-1}(q)$, even  if $P$ is non-compact.)\vspace {1mm} 

 {\it The   $m$-waist}  (mod 2), denoted   $waist_m (h)$, of  a homology class $h\in H_n(X;\mathbb Z_2)$ is\vspace {0.7mm}

 \hspace {30mm}  {\it the infimum of the  numbers $w$,} \vspace {0.7mm}

\hspace {-6mm}  such that $X$  receives a  Lipschitz  map from a compact  $m$-sliced cycle, $\phi: P^n\to X$, which represent  $h$, i.e. 
 $$\phi_\ast[P]=h$$ 
and  the

  \hspace {16mm}{\it the   images of all slices in $X$ have $m$-volumes $\leq w$}, \vspace {0.7mm}

   \hspace {-6mm} where these  "volumes  of the images"  are counted with multiplicities (which is unneeded for generically 
  1-1 maps.)
 \vspace {1mm}

 {\large \sf  \textbf {A$_{++}$}.} {\large{\sf Conjecture.}}  Let $X$ be a {\it closed} $n$-dimensional Riemannin manifold the scalar curvature of which is bounded from below as earlier:
 $$Sc(X)\geq n(n-1)(=Sc(S^n)).$$
  {\sl  Then  the slicing area of  the fundamental homology class $[X]\in H_n(X;\mathbb Z_2)$ is    bounded by
 $$waist_2[X]\leq const_{++}.$$}
(Ideally, one expects $$const_{++} = waist_2(S^n)$$
where $  waist_2[S^n]=area(S^2)=4\pi$ by an  Almgren's theorem.)

 \vspace {2mm}
 The above conjectures can be interpreted as  saying  that $X$ contains  "many" small  subsets of dimensions 1 and/or 2.

For instance,  {\large \sf A} implies that  that $X$ contains a {\it topologically significant/representative} family of
  $1$-{\it dimensional} subsets (graphs) with diameters  $\lessapprox\frac {1}{\sqrt\sigma}$.

  This   suggests the following.\vspace {1mm} 

 {\large \sf  \textbf {(a)} Conjecture}. {\sl If  $Sc(X )\geq \sigma>0$ and if $X$ is  a {\sf closed} (compact without boundary) manifold, then $X$ contains a closed minimal geodesic of length  $\leq \frac  {const_n}{\sqrt \sigma}$, or, at least, a stationary one-dimensional $\mathbb Z_2 $-current  of    diameter  (better length) $\leq \frac  {const_n}{\sqrt \sigma}$.}

\vspace {1mm}

And  {\large \sf   {A$_{++}$}}  actually  implies   the following.\vspace {1mm}

{\large \sf  \textbf {(a$_{++}$)} Conjecture}. {\sl Closed manifolds  $X$ with  $Sc(X )\geq \sigma>0$ contain  closed minimal surfaces  (i.e. stationary two-dimensional $\mathbb Z_2 $-currents) of areas 
$\leq \frac  {const_n}{\sigma}$.}\vspace {1mm}

Below is a weaker version of   {\large \sf {A}}  which already imposes  non-trivial topological constraints on $X$.
  \vspace {1mm}

  {\large \sf  \textbf {A$_-$}}.  {\large{\sf Conjecture.}} {\sl If $Sc(X)\geq n(n-1)$ then the {\sf filling radius} of
  $X$ is bounded by   
 $$ fil.rad (X)\leq const_-.$$}

 {\it Definition of {\sf fil.rad}.} If $X=(X,g)$ is {\it closed}    Riemannian  manifold then 
 the {\it filling radius} is equal to the infimum of  $R> 0$, such that the cylinder $X^\times =X\times [0,1)$ 
 admits a Riemannin metric $ g^\times $ with the following three properties.\vspace {1mm}
 
 $\bullet_1$ the restriction of $\hat g$ to $X= X_0\times \{0\}\subset X\times [0,1)=X^\times $ is equal to $g$; moreover, 
 $$ dist_{ g^\times }|X =dist_g.$$
This  means that the $g$-shortest curves in $X$  between all pairs of  points in $X$ minimise the  $ g^\times $-lengths of such curves  in $X^\times \supset X$.\vspace {1mm}

 $\bullet_2$  All points in $X^\times$ lie within distance  at most $R$ from $X$,  $$dist_{g^\times } (x^\times, X)\leq R\mbox{    for all }x^\times \in  X^\times .$$
 
 $\bullet_3$  The  $n$-dimensional volumes of the submanifolds $X\times \{t\} \subset  X\times [0,1)=X^\times$, $t<1$, with respect to $g^\times$ vanish in the limit for $t\to 0$,
 $$  vol (X\times \{t\})\to 0\mbox {  for } t\to 1.$$
 
  (The equivalence of this definition to  the usual one follows from  the the {\it filling volume  inequality}  see [9]\footnote {  [9]  L. Guth, {\sl  Notes on Gromov's systolic estimate}, Geom Dedicata (2006) 123:113-129. } and references therein).
 
 Then the filling radius  of a  compact manifold $X$ with boundary --  our  manifolds may, a priori, have boundaries and/or to be  incomplete -- is defined as $fil.rad$ of the double of $X$ along the boundary and $fil.rad$ of an open $X$ is defined via  exhaustions of $X$ by compact submanifolds.
 \vspace {1mm} 
 
 It is obvious that 
 \hspace {1mm}  {\large \sf  \textbf {A$_+$}} $\Rightarrow $ {\large \sf  \textbf {A}} $\Rightarrow $  {\large \sf  \textbf {A$_-$}} and that  {\large \sf  \textbf {A$_+$}} is optimal in a way.
 
Indeed, the product $X_r=X_0\times S^2(r))$, where $X_0$ is,  a  compact manifold and $S^2(r))$ is the 2-sphere  of small  radius $r\to 0$, (these spheres have 
$Sc(S^2(r))= \frac {2}{r^2}$),
has  $Sc(X_r)\geq (\frac {2}{r^2} - const_{X_0}) \to+\infty$, while the $(n-2)$-dimensional size/spread of $X_r$
 is  as large as that of $X_0$.

  Also one knows  (see   [17] at the end of this section and references therein)   that 
  $$\mbox{ {\large \sf  \textbf {A$_{++}$}} $\Rightarrow $ {\large \sf  \textbf {A$_-$}}}.$$
(It is plausible in view of    [18] that {\large \sf  \textbf {A$_{++}$}} $\Rightarrow $ {\large \sf  \textbf {A}}.)
    \vspace {1mm} 
 
 On the other hand, it is not hard to show that  if  the  {\sl if  the   isometry group of    a    Riemannin manifold  $\hat X$  acts  {\sf cocompactly} on  $\hat X$,
 i.e $\hat X/isom(\hat X)$ is compact, and if  $\hat X$ is  {\sf contractible},  then 
 $$fil.rad(\hat X)=\infty.$$}\vspace {1mm} 
 
Therefore, {\large \sf  \textbf {A$_-$}} yields  the following topological $Sc>0$-non-existence  corollary.\vspace {1mm}
 
{\large \sf  \textbf {B}. Conjecture.}  {\sl Closed manifolds $X$ with contractible universal coverings $\tilde X$ admit no metrics with $Sc>0$.}

(Granted {\sf \textbf { B}}, the non-strict inequality   $Sc(X)\geq 0$ implies that    $X$ {\it  Ricci flat} 
 by   Kazdan-Warner's perturbation theorem  (see [10] \footnote {{\it If a metric $g_0$ with $Sc\geq 0 $ {\it can't be perturbed} to $g$ with $Sc(g)>0$, then  $Ricci(g)=0$.}
 
[10]   J Kazdan, F. Warner,
  {\sf Existence and Conformal Deformation of Metrics With Prescribed Gaussian and
Scalar Curvatures,} 
Annals of Mathematics,  101,\# 2. (1975), pp. 317-331.}).   
And since $\tilde X$ is contractible, the universal covering   $\tilde X$ is {\it isometric to the Euclidean space $\mathbb R^n$,} $n=dim (X)$, by  the {\it Cheeger-Gromoll splitting theorem}.)
 
 \vspace {1mm}

\vspace {1mm}

 \hspace {30mm} {\large \sf \textbf {Remarks and  References.}}   \vspace {1mm}

However  plausible,  none of the    {\large \sf  \textbf {A}}-conjectures  (above dimension 2)  has been confirmed  except for 
  {\large \sf  \textbf {A$_+$}} for  $3$-manifolds $X$  with 
(apparently non-sharp) constant $const_+=2\pi \sqrt 6$ (see [14] below).

 On the other hand {\sf \textbf { B}} is known to hold  for many manifolds $X$,
starting from the case of $n$-tori due to Schoen and Yau.  
Later $B$ was proven  by a use of {\it  twisted Dirac operators}\footnote {This means: {\sl Dirac operators with coefficients in some (possibly infinite dimensional) vector bundles.}} 
 for several  classes of  manifolds with "large" universal coverings    
including those $X$ which admit metrics with non-positive sectional curvatures.

Below are a few  relevant papers  where one can find further references.

\vspace {1mm}

[11]  Yau, S.T., and Schoen, R.  {\sl On the Structure of Manifolds with positive Scalar Curvature}. Manuscripta mathematica 28 (1979): 159-184.\vspace {1mm}

 [12] J. Lohkamp,
{\sl The Higher Dimensional Positive Mass Theorem II}, (2016)
 	arXiv:1612.07505.\vspace {1mm}

[13] R. Schoen, S.T. Yau, {\sl Positive Scalar Curvature and Minimal Hypersurface Singularities},
 (2017) 	arXiv:1704.05490. \vspace {1mm}

In [11], the authors introduced their method of {\sl induction descent by minimal  hypersurfaces} and proved non-existence of metrics with $Sc>0$ on the $n$-tori\footnote {This trivially implies non-existence of  compactly supported perturbations with $Sc>0$  of the  Euclidean metric on $\mathbb R^n$.} and, more generally, on   $n$-dimensional manifolds $X$
which admit smooth maps $X\to \mathbb T^{n-2}$, such that the   homology classes in $H_2(X)$ represented by  the 
pull backs of  generic points  are {\it non-spherical.} 

Originally, this  method was limited to $n\leq 7$, but the techniques developed in [12] and [13]  apparently remove this limitation.
\vspace {1mm}

[14] M. Gromov, H. B. Lawson, Jr., {\sl Positive scalar curvature and the Dirac operator on complete Riemannian manifolds},  Publ. Math. IH\'ES 58 (1983), 295-408. 

In this paper besides  above mentioned  {\large \sf  \textbf {A$_+$}} for  $3$-manifolds,  
  we   {\it rule} out  complete metrics with $Sc>0$  on certain classes of manifolds, including   \vspace {1mm}

  {\sl  closed orientable $n$-dimensional spin\footnote{A manifold of dimension $n\geq 3$ is spin if the restrictions of the tangent bundle $T(X)$  to all immersed surfaces in $X$ are trivial bundles.

Most (all?) known  non-existence results for $Sc>0$ obtained for {\it spin manifolds} more or less automatically  generalise to manifolds whose {\it universal coverings are spin}, i.e where $T(X)$ trivialises on all immersed $2$-spheres in $X$.}
manifolds  $X$ which admit continuous  maps to complete  manifolds $Y$ with non-positive sectional curvatures,  such that the fundamental classes $[X]\in H_n(X)$ go to  {\it non-zero } classes in $H_n(Y)$ under these maps.}
\vspace {1mm}

[15] M. Gromov, {\sl Positive curvature, macroscopic dimension, spectral gaps and higher signatures,}  in Proc of 1993 Conf. in Honor of of the Eightieth Birthday of I. M. Gelfand, 
 Functional Analysis on the Eve of the 21st Century: Volume I  Progress in Mathematics,  (1996) pp 1-213, Vol. 132, 
 
This paper presents a   geometric  perspective on the Dirac operator and soap bubble methods  in the 
study of scalar curvature and related problems.
 \vspace {1mm}

[16] S. Markvorsen, M. Min-Oo,  {\sl Global Riemannian Geometry: Curvature and Topology}, 2012
Birkh\"auser.
 
 A chapter in this book \footnote {Also see MinOo,  {\sl K-Area, mass and asymptotic geometry}, 
 
 http://ms.mcmaster.ca/minoo/mypapers/crm }  offers a friendly introduction to the Dirac operator methods in the $Sc> 0$ problems.
 
\vspace {1mm}

 [17] L. Guth, {\sl Metaphors in systolic geometry}. In: Proceedings of the International Congress of Mathematicians.  2010,  Volume II, pp. 745?768.

[18] L. Guth,{\sl  Volumes of balls in Riemannian manifolds and Uryson width}.
Journal of Topology and Analysis,
Vol. 09, No. 02, pp. 195-219 (2017).

These two  papers and references therein give a fair idea of results  and ideas around the filling radius.  \vspace {1mm}

[19] D.Bolotov, A. Dranishnikov.  {\sl On Gromov's conjecture for totally non-spin manifolds}, (2015)
arXiv:1402.4510v6.

[20] M. Marcinkowski, {\sl Gromov positive scalar curvature conjecture and rationally inessential macroscopically large manifolds}.
 Journal of topology  9, 1; 105-116  (2016).
 Oxford University Press

The authors of these two papers are concerned with  topological versions of  {\sf \textbf { A$_+$}} for certain classes of  manifolds $X$.
\vspace {1mm}

[21] J. Rosenberg,  {\sl Manifolds of positive scalar curvature: a progress report}, in: Surveys on Differential Geometry,
vol. XI: Metric and Comparison Geometry, International Press 2007.

This is  survey of      topological obstructions  to metrics with  $Sc>0$ on spin manifolds  $X$ 
expressed in terms of   indices of   {\it Dirac operators twisted with $C^\ast$-algebras} of $\pi_1(X)$.

 Also 
obstructions for  {\it $4$-dimensional} manifolds $X$  with   non-vanishing {\it Seiberg-Witten invariants} due to  Taubes and Le Brun are described in this paper.

[22] M. Gromov, {\sl Morse Spectra, Homology Measures, Spaces of Cycles and Parametric Packing Problems.}
\url{www.ihes.fr/~gromov/PDF/Morse-Spectra-April16-2015-.pdf}

This is  an overview of waists and related invariants which may  bear some relevance  to $Sc\geq \sigma$. 

\vspace {1mm}

 \section {Extremality and  Rigidity  with Positive Scalar Curvature.}

 The proof(s) of  the above  { \sf {A}}-conjectures (let them be only  approximately true)  would  require  {\it constructions} of certain maps   or spaces  which makes these conjectures   difficult. 
 
What is  easier is getting  upper  bounds on the    "size" of an $X$ with $Sc(X)\geq \sigma>0$   by proving {\it lower bounds} on   dilations of {\it topologically significant} maps from $X$ to  (more or less) standard manifolds $Y$.

 The first {\it sharp} bound of this kind was proved in \vspace {1mm}
 
 \hspace {2mm} [23]  M. Llarull,   {{\sl Sharp estimates and the Dirac Operator},

\hspace {2mm} Math. Ann.
310 (1998), 55-71,}

\hspace {-6mm}followed by  \vspace {1mm}

 \hspace {2mm} [24]  M. Min-Oo, {\sl Scalar Curvature Rigidity of Certain Symmetric Spaces}, 
 
 \hspace {2mm} Geometry, Topology and Dynamics
(Montreal, PQ, 1995), CRM Proc. 
 
  \hspace {2mm} Lecture Notes, 15, Amer. Math. Soc., Providence, RI, 1998, pp. 127-136.

\hspace {-6mm} and  by\vspace {1mm}
 
  \hspace {2mm}  [25]  S. Goette and U. Semmelmann, {\sl Scalar curvature estimates for 
 
   \hspace {2mm}  compact symmetric spaces}, 
 Differential Geom. Appl. 16(1), (2002) 65-78,

\hspace {-6mm} where further references can be found.\vspace {1mm}

What is proven in these papers can be expressed in the  the following terms. \vspace {1mm}


  {\large\sf Extremality/Rigidity.}  A Riemannian metric  $\underline g$
  on a   manifold $Y$ is called {\it length extremal}  if it 
  
  {\sl can't be enlarged without making the scalar curvature smaller somewhere.}  
  \hspace {-6mm} Namely,
   the inequalities 
  $$ Sc(g)\geq Sc( g_0) \mbox { and } g\geq g_0$$ 
   for a  Riemannian metric $g$ on $Y$
   imply 
     $$Sc(g)=Sc(\underline g).$$
Then  the stronger implication
$$[Sc(g)\geq Sc(g_0)] \& [g\geq g_0]\Rightarrow [g =\underline g]$$ 
is qualified as  {\it length rigidity} of $\underline g$.\footnote{ Extremal manifolds define, in a way, the boundary of the domain $\{{\cal SC}\geq 0\}$  of manifolds with $Sc\geq 0$.}
\vspace {1mm}

{\it CY-Example.} If a closed manifold $Y$ admits no metric with $Sc>0$, then all  $g_0$ with $Sc(g_0)=0$ \footnote {The condition $Sc(g_0)=0$ implies  $g_0$ $Ricci(g_0) =0$ on these $Y$  by the Kazdan-Warner perturbation theorem, see [10] in section 3}  are extremal according to this definition.

Instances of such {\it scalar flat} manifolds    are flat  Riemannin manifolds (with universal coverings  $\mathbb R^n$)  and also (simply connected) 
 hypersurfaces $Z\subset \mathbb CP^{n+1} $  of  degree   $n + 2$ and even $n$, with {\it  Ricci flat   Calabi-Yau metrics,}
 where  non-existence of metrics  with $Sc>0$ on these $Z$ follows from the Lichnerowicz, theorem.

\vspace {1mm}

Next, define {\it area extremality} and {\it area rigidity}  
by relaxing the inequality 
$g\geq g_0$, which says in effect that 
 $ lenght_g(C) \geq  lenght_{\underline g} f(C)$
for all smooth curves $C\subset Y$), to
$$ area_g(\Sigma ) \geq  area_{\underline g} (f(\Sigma))$$
for all smooth surfaces  $\Sigma \subset Y$,   
 where the extremality and rigidity  requirements remains the same:
   $Sc(g)=Sc(\underline g)$ and $g =\underline g$.

   \vspace {1mm}

Stronger versions  of these extremalities and rigidities  allow modifications  of the topology as well as geometry of $Y$, where the role of  "topologically modified"  $Y$  are played  by a Riemannin manifold $X=(X,g)$ and    a map $f:X\to Y$, where the above inequalities are understood as 
$$ Sc(g)(x) \geq Sc (\underline g)(f(x)),   lenght_g(C) \geq  lenght_{\underline g} g(f(C))$$ 
\vspace {-1.6mm}  and \vspace {-1.6mm}
$$area_g(\Sigma ) \geq  area_{\underline g} (f(\Sigma))$$
correspondingly.
 
 Accordingly, the required conclusion for {\it extremality} is  
 $$Sc(g)(x) = Sc (\underline g)(f(x)),$$
while both, the  {\it length} and the {\it area rigidities}, signify that 
$$lenght_g(C) =  lenght_{\underline g} (f(C)).$$
for all smooth curves $C\subset X$.

 Of course, these definitions makes sense only for particular topological classes of manifolds $X$ and maps $f$, such for instance as the class $\{ {\cal DEG}\neq 0 \}$ of orientable manifolds of dimension $n=dim(Y)$   and {\it $C^2$-smooth} maps with {\it non-zero degrees}.  \vspace {1mm}
 
{\large \sf  \textbf { C}.  Problem.} {\sf Find verifiable  criteria for extremality and rigidity, decide which manifolds admit extremal/rigid  metrics and describe particular classes of extremal/rigid manifolds.}

 \vspace {0mm}

For instance,  \vspace {1mm}

{\sf do all closed manifolds which admits metrics with $Sc\geq 0$ also admit (length) 

extremal metrics?} \vspace {1mm}

More specifically, prove (disprove?) the following. \vspace {1mm}

  \vspace {1mm}

 {\large\sf \textbf { C$_1$}.  Conjecture.} {\sl All compact Riemannin symmetric spaces  are area extremal in the class  $\{ {\cal DEG}\neq 0 \}$ and those which have $Ricci>0$ (this is equivalent to absence of local $\mathbb R$ factors, and to  is finiteness of  fundamental group) are area rigid in this class.}  \vspace {1mm}

This  conjecture was proved by  Llarull  (see    [23]  above) in the case $Y=S^n$, under the additional assumption of $X$ being {\it spin}.\footnote {Since $\pi_1(SO(n))=
\mathbb Z_2$ for $n\geq 3$, there are at most two   isomorphism classes of vector bundles with $ rank \geq 3$ over connected surfaces $\Sigma$ (exactly two for closed $\Sigma$),  where the trivial bundle is called {\it spin} and where  bundles of $rank <3$ are    spin if their Whitney sums with trivial bundles are spin. An orientable  vector bundle  $V$ of   over a topological space $B$ is {\it spin} if the pullbacks of $V$ under continuous maps $\phi :\Sigma\to B$ for all surfaces $\Sigma$ are spin. A manifold $X$ is spin if its tangent bundle is spin. 

 The  spin condition is necessary for the definition of the Dirac operator on $X$ but some {\it twisted Dirac operators} make sense  on non-spin manifolds.}\vspace{1mm}

Then Min-Oo [24] proved {\sl area extremality}  for {\it Hermitian symmetric spaces} in the class  

\hspace {-6mm}  $\{{\cal SPIN},  {\cal DEG}\neq 0 \}$, where the maps $f:X\to Y$, besides  having  degrees $\neq 0$, must be {\it spin}.\footnote {A  map $f:X\to Y$   is {\it spin} if  the pullbacks  $\phi^!( T(X))$  for  maps of  surfaces,  $\phi:\Sigma \to X$, satisfy 

\hspace {27mm} [$\phi^!( T(X))$ is spin] $\Leftrightarrow$  $[(\phi \circ f)^!( T(Y))$ is spin]

\hspace {-4mm}for all 
$\Sigma$ and $f$.  Equivalently,  a map $f$ between orientable manifolds is spin  if the Whitney sum $T(X)\oplus f^!(T(Y))$ is spin.

Obviously, the identity map $id:Y\to Y$  is spin and  if $Y$ is spin, e.g. $Y=S^n$, then

 \hspace {37mm}     [$f:X\to Y$ is spin] $\Leftrightarrow [X$ is spin].   }

 This was generalised by Goette and Semmelmann [25] who proved \vspace {1mm}

 {\it area extremality in  $\{{\cal SPIN},  {\cal DEG}\neq 0 \}$ of  compact (here it means {\sl closed})  K\"ahler manifolds  with $Ricci \geq 0$, rigidity for $Ricci \geq 0$. } \vspace {1mm}
 
 Moreover, they establish \vspace {1mm}

  {\sl area rigidity in  $\{{\cal SPIN},  {\cal DEG}\neq 0 \}$ of certain (non-Hermitian)   compact symmetric spaces including those with non-vanishing Euler characteristics and also of Riemannian metrics on $S^{2m}$   with positive curvature operators.}
  \vspace {1mm}

 These extremality and rigidity theorems are    proven in the non-K\"ahlerian  cases  by   {\it sharply evaluating}   the  contribution from  $f^!(\mathbb S^+(Y))$ in the  Schroedinger-Lichnerowicz formula  for the Dirac operator on $X$ twisted with the $f$-pullback of  the {\it spinor$^+$ bundle} $\mathbb S^+(Y)$ which is, in the case where $\chi(Y)\neq 0$ is 
confronted with the index theorem.

(The  case of odd dimensional spheres  $S^{2m-1}$, which  depends on an additional argument(s)  applied to maps $X\times S^1\to S^{2m}$ \footnote{Llarull uses the product  metric on $X\times S^1$, where his calculation applies even though the scalar curvature   $Sc(X\times S^1)$, which is $\geq  Sc(S^{2m-1})$, may be smaller than   $Sc(S^{2m})$. 

Alternatively, one can use the {\it spherical suspension metric} $g^S$ (of $g$ on $X$)  on (the bulk of)  $X\times S^1$,   which has $Sc(g^S)\geq  Sc(S^{2m})$ and thus  allows a formal reduction of the $2m-1$ case to that of $2m$.}
seems to apply only  to metrics on $S^{2m-1}$  with constant sectional curvatures.)   
 
 And in the K\"ahler  case, this is done with the {\it "virtual square root" of the canonical} (complex) {\it line bundle} on $Y$ instead of  $\mathbb S^+(Y)$.
 
\vspace{1mm}

{\large \sf Spin or non-Spin?} In all of the above cases one can replace
 the spin condition for $f:X\to Y$ 
 by this condition for 
 the corresponding map between the  universal coverings,  $\tilde f: \tilde X\to \tilde Y$, where   a version of Atiyah's $L_2$-index theorem applies. 
 
 \vspace {1mm}

Probably,  

\hspace {6mm} {\large \sf  "spin" can be removed  all together in these theorems}  

\hspace {-6mm} but this  seems beyond reach of the present day methods.\footnote{Apparently, no single case of extremality of a  closed simply connected  manifold $X$ of dimension $n
\geq 3$ is amenable to the  the minimal hypersurface techniques, except, may be(?)  
for $X=S^3$.}

\vspace {1mm}

On the other hand, the spin condition is essential for the extremality in the class ${\{\cal SPIN}, {\cal DEG}_{\hat A}\neq 0\}$ where the dimension of $X$ can be greater than $n=dim(Y)$ and where
the condition $deg(f)\neq 0$ is replaced by $deg_{\hat A}(f)\neq 0$, where the {\it $\hat A$-degree}  $deg_{\hat A}(f)$ stands for 
 the $\hat A$-genus of the $f$-pull back of a generic  point $y\in Y$,
 $$deg_{\hat A}(f)=\hat A(f^{-1}(y)).$$ 
(Here, strictly speaking,  $f$ must be smooth;  if $f$  is just continuous, this applies to a smooth approximation of $f$, where the so defined  $\hat A$-degree  does not depend on a choice of approximation.).
 
 This implies for instance, that \vspace {1mm}

 {\sl the products of the above $Y$, e.g. of $Y=S^n$ by the Calabi-Yau manifolds with $\hat A\neq 0$, e.g with $Z$ from the above CY-example are area extremal in the class ${\{\cal SPIN}, {\cal DEG}_{\hat A}\neq 0\}$ as well s  in the class $\{\widetilde {\cal SPIN}, {\cal DEG}_{\hat A}\neq 0\}$ where spin condition is delegated to $\tilde f :\tilde X\to \tilde Y$}. \vspace {1mm}

\vspace {1mm}

Notice, however, that {\it neither} simply connected  Calabi-Yau manifolds $Z$ themselves   {\it nor} their
 products by $Y$ are extremal in the class ${\{\cal SPIN}, {\cal DEG}\neq 0\}$,  at least if $dim(Z)\geq 5$.
 
 Indeed the connected sums $X=Z\#(- Z)$, where "$-$" stands for  the reversal of orientation  and where the obvious map   $Z\#(- Z)\to Z$ has degree 1,   admit metrics with $Sc>0$ by Stolz'  theorem mentioned    in  section 2.
 \vspace {1mm}sleeker

 It seems that the there are two divergent, yet interconnected by bridges,   branches in the   tree of $Sc(X)\geq 0$, where  a smoother and sleeker  one  involves  differential  structure  and depends on spin, while   the other  one 
is made of rougher staff  such as the homotopy classes of $X$. \footnote {The smooth branch is manifested by $\hat A$ and the $\mod 2$ $\alpha$-invariant in the index formula while the rough branch is  
 represented by   the Chern character  and   supported by  minimal hypersurfaces. }

Probably,  the  second branch can be transplanted to a  harsh world inhabited by  singular spaces    
 but fully cleaning  off  spin  from this branch  is by no means  easy even for smooth $X$.

  \vspace {1mm}

 {\it  Extremality and Rigidity of Products.}  It seems not hard to show\footnote {I have not verified  the proof in detail at this point.} that the Riemannian products of the area extremal/rigid  manifolds in the above examples are  area extremal/rigid which suggests to the following.  \vspace {1mm}
 
  {\large\sf \textbf { C$_2$}. Question.} {\sf Are the  Riemannian products of all  area extremal/rigid  manifolds  area extremal/rigid}
 
  \vspace {1mm}
 {\it Smoothing Lipschitz Maps.}  The   length  extremal/rigid manifolds in some {\it homotopy} class of smooth maps remain 
 extremal/rigid  in the corresponding class of  Lipschitz maps  $f$.  
 
 This  can be  proven by a smooth  approximation of these $f$ with a minor change of their length dilations. 
 
 But  
 
 \hspace {3mm} {\large \sf this is  unclear for the {\it area} extremality and/or  area rigidity}, 
 
 \hspace {-6mm}since, conceivably(?)  all smooth approximation  $f'$ of a Lipschitz map $f:X\to Y$  may have $area (f'(\Sigma))>> area (f(\Sigma))$ for some $\Sigma$.\vspace {1mm}

 {\it Normalisation by Scalar Curvature: Extremality/Sc and Rigidity/Sc}.    A map $f:X\to Y$ between Riemannin manifolds  $X=(X,g)$ and  $Y=(X,g_0)$  with positive scalar curvatures,  $Sc(g), Sc(g_0)>0,$ is called {\it length decreasing/Sc} if it decreases the length of the curves measured in the metrics $Sc(X)^{-1}g $ and $Sc(g_0)^{-1}g_0$, i.e. if it decreases the integrals of $\sqrt {Sc}$ over all curves in $X$.
Similarly one understand  {\it decrease/Sc of areas} of surfaces $\Sigma\subset X$ under maps $X\to Y$, etc.\footnote{It may (or may not) be worthwhile to normalise by  $g \leadsto  n(n-1) Sc(X)^{-1}g$, $n=dim(X)$, and see what happens for $n\to \infty$.}

 
Accordingly, one defines {\it length/area extremality/Sc} of a $Y$ as non existence of strictly length/area  decreasing/Sc maps $X\to Y$ in a given class of manifolds and maps, while the 
rigidity/Sc signifies that all length/area non-increasing/Sc maps $f:X\to Y$ are homotheties (similarities)   with respect to the original metrics, i.e. $f^\ast(g_0) =const\cdot g$.

 Since the   "contribution of the twist" to the  Schroedinger-Lichnerowicz formula for the twisted Dirac opertor  on $X$  scales  as 
 $Sc(X)^{-1}$, the arguments from the above cited papers  based on this formula deliver
 the corresponding extremality/Sc and rigidity/Sc results. (This was pointed out  in [26]\footnote{[26]  M.Listing, {\sl   The  Scalar curvature on compact symmetric spaces},
 arXiv:1007.1832, 2010 - arxiv.org.})

 {\it Category ${\cal R}_+/_{Sc}$}.  Let this be the category of Riemannian manifolds with $Sc> 0$ 
 and length (alternatively, area) non-increasing/Sc maps. \vspace {1mm}
 
  {\large \sf \textbf {C}$_3$.  Question.} {\sl How much of the geometry of spaces with $Sc>0$ can be  
 reconstructed in the category theoretic  language of   ${\cal R}_+/_{Sc}$?} \vspace {1mm}

 \vspace {1mm}

{\it Extremality beyond $Sc\geq 0$}. The condition $Sc(g)\geq 0$ may be  not indispensable for  extremality of $g$.  

For instance, the  double of the unit hyperbolic disk is (kind of) extremal for the natural $C^0$-continuous metric  on it 
and there are similar high dimensional examples. But it is unclear if such metrics are ever  smooth. 

 

\vspace {1mm}

 {\it Relativisation of Non-existence Theorems for  $Sc>0$.} Let $Y$ be a closed length or  area  {\it extremal or rigid} manifold  in some  class of smooth manifolds $X$ and smooth maps $f: X\to Y$,  
 where this class  is invariant under homotopies of maps.

Then,  most (all?)  known  { \it Dirac operator obstructions} to the existence of   metrics with $Sc>0$ on  closed  manifolds  $X_0$ {\large \sf naturally extend} to  {\it similar obstructions} to the existence of (strict)  {\it  area decreasing/Sc}  maps in  certain   homotopy invariant  classes of maps $X\to Y$, including $X=X_0\times Y\to  Y$ for $(x_0,y)\mapsto y$.
 
 \vspace {1mm}
 
 For instance, one knows that (co){\it homologically symplectic manifolds} $X_0$ with $\pi_2(X_0)=0$ admit no metrics with $Sc>0$ and the proof of this (see  [15] cited  in the previous section)    also implies that \vspace {1mm}words
 
{\sl   if   $Y$ is the above area-extremal manifold, e.g. $Y=S^n$, then no  {\it homologically symplectic}\footnote{A smooth proper map between  orientable manifold, $f: X\to Y$, is {\it homologically symplectic}   
 if  the difference of the dimensions $n_0=n-m$ for $n=dim(X)$ and $m=dim(Y)$  is {\it even}  and if   there exists a closed 2-form $\omega$ on $X$ such that the integrals of $\omega^\frac {n_0}{2}$ over the $f$-pullbacks of generic points $y\in Y$ do not vanish.
 
 In other words, 
 the real fundamental cohomology class  $[X]^\circ\in H_{comp}^n(X;\mathbb R)$  with compact support is equal to 
 the $\smile$product  of the $f$-pullback of $[Y]^\circ\in H_{comp}^m(Y, \mathbb R)$ and the $\frac {n_0}{2}$th$\smile$power of the class $[\omega]\in H^2(X;
 \mathbb R)$, 
 $$[X]^\circ = f^\ast ([Y]^\circ)\smile [\omega]^\frac {n_0}{2}. $$}
  map $f:X\to Y$, which, moreover, induces an  {\it isomorphism} $\pi_2(X)\to \pi_2(Y)$}, can be {\it strictly area decreasing/Sc}. 

 
This suggests the following.  \vspace {1mm}
\vspace {1mm}
 
 {\large \sf \textbf {C}$_4$.  Conjecture.}  Let $g$ be a metric on $X$ and $f_0: X\to Y$ be   a (smooth?) strictly length (area?) decreasing/Sc  map in this class.
 
 {\sf Then  there exists a smooth  map $f$  homotopic to $f_0$ transversal to  a point $y_0\in Y$,  such that  the $f$-pullback  submanifold  $f^{-1}(y_0)\in X$ admits a metric with $Sc>0$.}\vspace {1mm}
 
 Also other properties, e.g. extremality, of   {\it manifolds $X$} with $Sc(X)>0$ may have  counterparts for   length and area decreasing/Sc {\it maps $X\to Y$}  and, furthermore, for {\it foliations} on $X$. \vspace {1mm}

 {\large \sf \textbf {C}$_5$.  Question.} {\sl Are   infinite dimensional counterparts of compact  symmetric spaces, e.g.  the Hilbert sphere $S^\infty$, extremal/rigid in some class(es) of perturbations of their metrics?}

 \section {Extremality and Gap Extremality of Open manifolds.}

 Let $U\subset Y$ be an open subset in a extremal or rigid  Riemannin   manifold $Y$ where  the extremality/rigidity  for this $Y$ follows by the   twisted Dirac operator  argument  from the previous section. Then the same  argument yields the following. \vspace {1mm}
 
 {\Large  $\star$ }  {\sl If the complement $Z=Y\setminus U$ is non-empty, yet  {\sf LC-negligible}}  (explained below)   {\sl then   {\large \sf no complete} orientable  Riemannian manifold admits a smooth
area non-increasing/Sc  map $f:X\to U$, which has
 {\it non-zero degree}\footnote {Maps $f:X\to Y$ of non-zero degree, by definition,  must be equidimensional and proper.} and the lift of which to to the universal coverings, $\tilde f:\tilde X\to \tilde U$, is {\large \sf spin}.}

 \vspace {1mm}
 
 {\it LC-negligible Sets .}  A piecewise smooth polyhedral subset $Z$ in a Riemannin manifold $Y$ is called  LC-negligible if the Levi-Civita connection on the tangent bundle of $X$ restricted to $Z$ is {\it split trivial.}
 For instance,\vspace {1mm}

  $ \bullet $   finite subsets in $Y$  are LC-negligible;

 $ \bullet $ piecewise smooth  graphs  $Z\subset Y$ with trivial monodromies around the cycles, e.g. disjoint unions of trees, are   LC-negligible;

  $\bullet $ simply connected  {\it isotropic} (e.g.{\it  Lagrangian})  submanifolds in K\"ahler manifolds are LC-negligible.\vspace {1mm}

 This definition extends to   general closed subsets    $ Z$, such as {\it Cantor sets}, for instance,  by  requiring that the monodromies along smooth curves $C$  in the $\varepsilon$-neighbourhoods of $Y$ are  $o (\varepsilon \cdot$ length $(C))$  as $\varepsilon \to 0$ but the geometry behind this definition needs to be clarified.

 \vspace {1mm}

 {\large \sf \textbf {D}$_1$.  Problem.}  {\sl Study essential properties, such as {\it the Hausdorff dimensions},   of these subsets $Z\subset Y$ and  find cases (if there are any)    where {\Large  $\star$ } remains valid for 
 small, yet non-LC-negligible  $Z\subset Y$, e.g. for (generic) smooth  curves $Z$ in $Y$}.\vspace {1mm}
  
Notice in this regard  that a simple surgery type argument  (see Stolz' paper  [5] cited  in  section 2 and references therein) shows that 
\vspace {1mm} 

$\bullet$ if $Z$ is equal to the  $k$-skeleton ${\cal T}^k$  of a smooth triangulation $\cal T$   of a   compact Riemannin manifold $(Y, g_0)$, for  $k\geq 2$,  then  $U=Y\setminus Z$ admits a  complete metric  $g\geq g_0$ with  $Sc(g)\geq \sigma_0=\sigma_0(Y, Z)>0$.

Moreover, it is easy to show that \vspace {1mm}

{\sl the  complements $U_\varepsilon=Y\setminus {\cal T}_\varepsilon^k$  of the $k$-skeleta  of the "standard fat" $\varepsilon$-refinements\footnote{It is more practical  to start with a  {\it cubilation}  ${\cal T}$ of  $Y$ which can be  canonically $\varepsilon$-refined  for $\varepsilon =\frac {1}{i}$, $i=2,3,...,$ by subdividing each $m$-cube  into $i^m$-sub-cubes in an obvious way.} of $\cal T$  admit complete Riemannin metrics $g_\varepsilon\geq g$ the scalar curvatures of which  for  $k\geq  2$ satisfy
$$Sc(g_\varepsilon ) \geq const\frac {1}{\varepsilon^2}$$
for some constant $const=const (Y,{\cal T})>0.$}

Thus   {\Large  $\star$ } fails to be true,  for   $Z={\cal T}_\varepsilon^k$, $k\geq 2$, and small (how  small?)  
$\varepsilon$.
\vspace{1mm}

On the other hand,  the {\sl torical band width inequality}  from  the next section shows that   if, for instance,  $Z$ is a codimension two torus in $Y$, e.g. $Z=\mathbb T^2\subset S^4$, then
  the complement $U=Y\setminus Z$ admits no complete metrics with $Sc\geq \sigma >0$ whatsoever
  and the same applies to a large (how large) class of codimension two polyhedra  $Z\subset Y$ with contractible universal coverings. 
   
   \vspace{2mm}

 Non-existence   of  {\it complete} metrics $g\geq g_0$ with $Sc>\sigma_0$ on the above    $U= (U, g_0)$  with $Sc(g_0)=\sigma_0$ may be interesting in its  own right but this  can't be regarded as 
extremality of  $g_0$, since  a comparison of the manifolds   $(U, g_0)$, which  have {\it bounded diameters} with  their competitors   $(U,g)$ of {\it infinite size} 
is patently unfair.    The true  extremity issue for these $U$, thus,  remains    unresolved.\vspace {1mm}

 {\large \sf \textbf {D}$_2$.  Question.} Do there ever  exist   length extremal  domains $U\subset Y$, $U\neq Y$,  in closed connected  Riemannin manifolds $Y$ of dimensions $\geq 3$?

For instance, is   the  the sphere $S^3$ minus a point (or the $3$-torus minus a point) extremal?

\vspace {1mm}

We still do not know the answer but, on the other hand,   the following {\it warped product}  construction sometimes delivers examples 
of both complete and non-complete  extremal and rigid manifolds (compare \S12 in [14] cited in section 3 and [27] cited  below).

\vspace {1mm}

Let $Y_0 =( Y_0,g_0)$ be a Riemannian  manifold   with constant scalar curvature $\sigma_0$ and let 
$g_1=\varphi^2g_0 + dt^2 $  be a Riemannin metric on $Y_1=Y\times (l_-,l_+)$ for $-\infty \leq l_-<l_+\leq  \infty$,  for some  smooth function  $\varphi=\varphi(t)>0$ for $l_-< t <l_+ $. 

Then, by elementary calculation, 
$$ Sc(g)= \frac {\sigma_0}{\varphi^2} -2n\frac {\varphi''}{\varphi}-n(n-1)\frac {\varphi'^2}{\varphi^2},\mbox{ where } n=dim(Y_0). \leqno {\RIGHTcircle}$$

Now,  let $g$  have {\it constant  scalar curvature}, say $Sc(g_1)=\sigma_1$ for a given $\sigma_1\geq 0$, and prescribe:  $\varphi(0)=1$ and $\varphi'(0)=0$.

Then  ${\RIGHTcircle}$, regarded as an  ODE and  rewritten as 
$$f'' =-\frac {1}{2}(n+1)f'^2 +\frac {\sigma_0}{2ne^{2f}}-\frac {\sigma_1}{2n}\mbox { \hspace {1mm}   for } f =\log \varphi,$$
admits a unique solution $f$ on some maximal (extremal)  open interval $(l^{ext}_-,l^{ext}_+)$ beyond which the solution does not extend.\vspace{1mm} 

{\it Examples.}(a)  If $Y_0=S^n$  and  $\sigma_1=n(n+1)$, then $Y_1$ is equal to $S^{n+1}$ minus two opposite points.\vspace{1mm} 

(b) If $Y_0=\mathbb R^n$ and   $\sigma_1=0$, then $Y_1=\mathbb R^{n+1}$.\vspace{1mm} 

(c) If $Y_0=\mathbb R^n$, $\sigma_1=n(n+1) =Sc(S^{n+1})$  and  $n=1$, then $Y_1$ is equal the universal covering of   $S^2$ minus two opposite points. 

In general,  the manifold $(Y_1,g_1)$ is   uniquely characterised  by the following  three properties.

\vspace {1mm}

$[\Circle_{n(n+1)} ]$ The scalar curvature of $Y_1$  is everywhere equal to $n(n+1)$  for $n=dim(Y_1)-1$.  \vspace {1mm}

$ [\Circle_{ O(n)\rtimes \mathbb R^n }]$  The isometry group of   $Y_1$ is    $Iso(\mathbb R^n)=
 O(n)\rtimes \mathbb R^n$ times $\mathbb Z_2$. (This  $\mathbb Z_2$ corresponds to the involution $t\leftrightarrow -t$.) \vspace {1mm}

$[\Circle_{2\pi/n+1}]$  The  {\it band  width}  of  $Y_1$ is  $\frac{2\pi}{n+1}$, where this width  is understood in the present case    as the distance between the  two (one point)  boundary components of $Y_1 $ in the  metric completion $\bar Y_1\supset Y_1$. 

(The band-like shape of  $Y_1$  is 
best seen for $dim(Y_1)=2 $, where  this  $Y_1$ is equal to the  universal covering of the doubly punctured sphere  $S^2$.)

Alternatively, one might say that the in-radius of $Y_i$ is equal to $\frac{\pi}{n+1} $:\vspace {1mm}

\hspace {5mm} {\sl there are closed  {\it compact } balls in $ Y_1$ of all  radii $R<
\frac{\pi}{n+1} $ but no ball of 

\hspace {5mm} radius $\geq \frac{\pi}{n+1}$ is compact.}\vspace {1mm}


\vspace {1mm}

{\it Gap Extremality.} We do not know if  the above  spheres minus pairs of points  are extremal for $n\geq 2$ but the Euclidean spaces $\mathbb R^{m}$  are definitely {\it not length extremal} starting from  $m=2$. 

In fact, there are (obvious, $O(m)$-invariant) 
metrics $g\geq g_{Eucl}$ on $\mathbb R^{m}$ with $Sc(g_1)>0$ for all $m\geq 2$.

On the other hand, 
\vspace {1mm}

({\Large $\ast$})  {\sl no metric $g\geq  g_{Eucl}$ on $\mathbb R^m$ may have  $Sc(g)\geq \varepsilon>0$}.(See [15] cited in section 3.)
\vspace {1mm}

This suggests the following weaker version of extremality for non-compact manifolds which we call 
{\it gap extremality.} 

A metric $g_0$ on $Y$ is {\it  $\varepsilon$-gap  length extremal} if {\it no} $g\geq g_0$ on $Y$  satisfies 
 $$Sc(g)- Sc(g_0) > \varepsilon.$$
Then  $g_0$ is called {\it gap  length  extremal} if it is   $\varepsilon$-gap length extremal  for all $\varepsilon>0$ ($0$-gap extremal=extremal).

Similarly one defines {\it area gap extremality} and {\it gap extremality for classes of maps} $f:X\to Y$.
(But I am not certain what  a {\it workable} definition of normalised  gap extremality/Sc should be.)  

\vspace {1mm}
Whenever the  twisted Dirac operator  argument  from the previous section yields area extremality of a closed manifold  $Y$, e.g. if $Y=S^n$ or $Y=\mathbb CP^n$,  this argument, combined with that from [15] (cited in section 3) for $\mathbb R^m$,  also delivers \vspace {1mm}

({\Large $\ast\ast$}) {\sl gap area  extremality of $Y_m=Y\times \mathbb R^m$ for all $m=1,2,...$, as well as this extremality for  smooth proper spin maps $f:X\to Y_m$ of non-zero degrees.}  \vspace {1mm}

{\sl if a smooth proper spin map $f:X\to Y_m$ of non-zero degree  decreases the areas of all surfaces $\Sigma\subset X$, 
 then, given $\varepsilon>0$,  there exists a point $x\in X$, such that 
 $$Sc(X)(x)- Sc(Y')(f(x))< \varepsilon.$$}\vspace {1mm}

 {\large \sf \textbf {D}$_3$. Question.} Does gap extremality is always stable under  $Y\leadsto Y\times \mathbb R^m$? (Beware of $dim(Y)=4$.)\vspace {1mm}

One can't discard of $\varepsilon$ for $m\geq 2$  but  the true area (or, at least length) extremality of $Y'=Y\times \mathbb R $  (that allows $\varepsilon=0$)  may be provable by some  twisted  Dirac operator argument. 
For instance,  if $Y=\mathbb T^n$ this follows from theorem 6.12 in   [14] (cited in section 3). 
Alternatively, one might  use  minimal hypersurfaces and soap bubble in $X$ the $f$-images of which separate  the two ends in $Y'=Y\times \mathbb R$ but then onr  would   face a  possibility of non-compact minimal hyper surfaces in $X$  and  would be obliged to resort to  imposing extra assumptions on $X$,  e.g.  uniform two sided bounds 
on the sectional curvatures of $X$.\vspace {1mm}

Finally, let us  look at  the    manifold $Y_1$, which has  the band width $\frac {2\pi}{n+1}$, in the above Example  (c).   

It is plausible  that this  $Y_1$ is length gap extremal but not length  extremal starting from $D=dim(Y_1)=3.$

And  what  we definitely know  is that  \vspace {1mm}

{\it  the quotient space $Y_1/\mathbb Z^n=\mathbb T^n\times (-\frac {\pi}{n+1}, \frac {\pi}{n+1})$, $n+1=dim(Y_1)$, is length

 extremal}.\vspace {1mm}

We shall see the reason for  this in the next section, where we  shall also explain the current  status of the rigidity problem for  these manifolds.

\section{ Bounds on Widths of Bands with Positive Curvatures.} 

Let us start with the following question which, on the surface of things, has nothing to do with scalar curvature.

Given a smooth $n$-dimensional manifold $X$ immersed\footnote{ A smooth map $X\to Y$ is an immersion if it is a diffeomorphism of  small neighbourhoods in $X$ to smooth submanifolds in $Y$.}
 into a complete Riemannian manifold $Y$
denote by   $rad^{\perp}(X\hookrightarrow Y)$ the maximal $R$, such  that  the normal exponential map  
 $$ exp^\perp: T^\perp(X)=T(Y)_{|X}\ominus T(X)\to Y,$$ 
 is locally injective on the subbundle $B^\perp(R)(X)\subset T^\perp(X)$ of open normal $R$-balls $B_x^{N-n}(R)\subset  T^\perp(X)$, $x\in X$.
 
  (If the ambient space  $Y=\mathbb R^n$, then  $rad^{\perp}(X;\mathbb R^n)$  is equal to the reciprocal of the supremum  of    the principal curvatures of  $X$.)  

Take the supremum of these  radii over  all immersions $f:X\hookrightarrow  Y$,
set
$$suprad^{\perp}(X;Y)=  \sup_f  rad^{\perp}(X\underset{ f}\hookrightarrow Y)  $$
and let 
$$suprad_N^{\perp}(X)=\sup_{f_\circ } rad^{\perp}(X\underset{ f_\circ}\hookrightarrow \mathbb R^N),$$  
where the latter "$\sup$" is taken over all immersion $f_\circ$  from $X$ to the {\it unit ball  $B^N(1)\subset \mathbb R^N$.}


(The notation $suprad^{\perp}(X; B^N(1))$ would be  unjustified, since the image of the exponential map 
may be not contained in   $B^N(1)$.)  \vspace {1mm}

 {\large \sf \textbf {E}$_1$. Problem.} Evaluate  $suprad_N^{\perp}(X)$ in terms of the topology of $X$.\vspace {1mm}

{\it Examples.} (a)  It is obvious that $suprad_N^{\perp}(X)\leq 1$ for all closed manifolds $X$, where the equality holds if and only if $X$ is diffeomorphic to $S^n$ and $N>n$.

(b) Let  $X_k$ is diffeomorphic to the product of $k$ spheres,
$$X_k= S^{n_1}\times ...\times S^{n_k}, \hspace {0.5mm}  n_k\geq 1.$$
Then $$suprad_N^{\perp}(X)\geq\frac {1}{ \sqrt k}\mbox  { for all } N\geq (n_1+1)+....+(n_k+1).$$

But we do not know, for instance,  whether
$$suprad_N^{\perp}(X_k)\to 0\mbox{   for  $N=dim(X_k)+1$ and $ k\to \infty. $}$$
or, on the contrary, if
$$suprad_N^{\perp}(X)\geq \rho_0 $$
for {\it all}  manifolds $X$, (e.g. for all $X_k$) all sufficiently large $N\geq N(X)$ and some universal constant $\rho_0>0$, say  $\rho_0=0.001$.

All known upper  bounds on  $suprad_N^{\perp}(X)$  -- {\sl am I  missing  something obvious}?     {\sf 

\hspace {10mm}exclusively apply to manifolds $X$ which admit no metrics with} $Sc>0$.

A simple way to obtain such a bound is  as follows.

1.Scale     $B^N(1)\to B^N(\frac {1}{2})$, project  $B^N(\frac {1}{2}) $  to $S^N$ from the south pole  of $S^N$ and observe that  this distorts the curvatures of submanifolds  $X$ in the ball  $B^N(1)$ by  a finite amount independent of $X$ and $N$.

2. Apply  the Gauss formula to $X \hookrightarrow  S^N$ and thus show that the supremum of the principal curvatures of $X$ in $S^N$ satisfies   
$$ supcurv(X\hookrightarrow S^N) \geq \frac{\sqrt{n-1}} {N-n}$$ and  therefore, 
 $$suprad_N^{\perp}(X)\leq const \cdot \frac {N-n}{\sqrt{n-1}}$$
{\sl for all $n$-dimensional manifolds $X$ which admit no metrics with $Sc>0$  and for   some constant $const \leq 100.$} (See  [27]  cited below for details.)

It follows, for instance,  that  there are  \vspace {1mm}

{\sl exotic spheres $\Sigma^n$ of   dimensions $n=9, 17,  25, 33,...$,
such that 
$$suprad^\perp _{n+1}(\Sigma^n)\leq\frac {100}{\sqrt {n-1}},$$}  \vspace {1mm}
but one has no idea how sharp this inequality is and  if  there are  similar inequalities  for   exotic spheres which admit metrics with $Sc>0$.\vspace {1mm}

The above also applies to tori $\mathbb T^n$, since these admit no metrics with $Sc >0$ either, but here the following  better (but, probably, still very  far from being sharp)  inequality is available.
$$suprad^\perp _{n+1}(\mathbb T^n)\leq \frac{ 2\pi}{n+1}.$$

This is proven again by passing to $S^{n+1}$, where all we use of the geometry of $S^{n+1}$ is the inequality $Sc(S^{n+1})\geq n(n+1).$
 (Isn't it  amazing that  there is no apparent {\it direct} proof of  a  {\it much  stronger} bound on 
 $rad^\perp (\mathbb T^n\subset B^{n+1}(1)$.))

Namely, the above  bound on  $suprad^\perp _{n+1}(\mathbb T^n)$ trivially  follows from the following.
\vspace {1mm}

{\large \sf Torical Band Width Inequality.} {\sl Let $g$ be a metric with $Sc(g)\geq n(n+1)=Sc(S^{n+1})$ on the torical band (cylinder)  $\mathbb T^n\times [-1,1]$. Then the distance between the two boundary components of this band satisfies
$$ dist_g(\mathbb T^n\times \{-1\}, \mathbb T^n\times \{1\}) <\frac  {2\pi}{n+1}.\leqno{  \mbox {
{\Large $[$}{$\varocircle$}{\small$ _\pm <  \frac { 2\pi}{ n+1}$}{\Large$]$}}} $$}
 
This is proven in [27]\footnote{[27] M.Gromov,
{\sl Metric Inequalities with Scalar Curvature.} 

\url {http://www.ihes.fr/~gromov/PDF/Inequalities-July\%202017.pdf}.} 
with  a relative version of the  Schoen-Yau minimal hypersurface method.

Besides a bound on  $suprad^\perp _{n+1}(\mathbb T^n)$, the inequality   {\Large $[$}{$\varocircle$}{\small$ _\pm <  \frac { 2\pi}{ n+1}$}{\Large$]$}  (trivially)  implies that
\vspace {1mm}

{\sl the warped product  metric $\varphi^2(t) g_{\mathbb T^n}+dt^2$ on  $\mathbb T^n\times (-\frac {\pi}{n+1}, \frac {\pi}{n+1})$  
with $Sc=n(n+1)$, which was introduced  in the previous section, is   length extremal.}\vspace {1mm}

Also, the  argument in  [27]  yields {\it length rigidity} of this metric for $n\leq 6$, while  the general case needs an  elaboration  on recent   results  on "irrelevance of singularities" of minimal hypersurfaces proved in the papers [12] and/or [13] cited in section 3.

\section {Extremality and Rigidity  of Convex  Polyhedra.}

Let   $P\subset  \mathbb R^n$ be a compact convex  polyhedron with non-empty interior, let   $Q_i\subset P$, $i\in I$, denote its $(n-1)$-faces and let 
$$\angle_{ij}(P)=\angle(Q_i,Q_j)$$
denote its dihedral angles.

Say that $P$ is {\it extremal}  if all  convex polyhedra $P'$ which are combinatorially equivalent to $P$ and which have 
$$\angle_{ij}(P')\leq \angle_{ij}(P)\mbox { for all }  i,j\in I,$$ 
necessarily  satisfy
$$\angle_{ij}(P')= \angle_{ij}(P).$$

 It is known  -- the proof is elementary -- that \vspace {1mm} 
 
  {\sl the simplices and the rectangular solids are extremal and also   all $P$ with
  
   $ \angle_{ij}(P)\leq \frac {\pi}{2}$,  are   extremal.}\vspace {1mm}

But  it is unclear (at least to the present author)  what are (if any)  non-extremal $P$.

What we are truly interested in, however, is  extremality (and rigidity) of $P$ under transformations which keep 
the faces $Q_i$  convex (rather than flat) or, even better,  {\it mean convex}, i.e. keeping  their mean curvatures non-negative.

Thus, we say that $P$ is {\it mean convexly extremal} if {\it there is no $P'\subset \mathbb R^n $ diffeomorphic to $P$}  and 
such that

$\bullet $  the  faces $Q'_i\subset P'$ corresponding to all $Q_i\subset P$ have $mean.curv(Q'_i)\geq 0$,

$\bullet $ the dihedral angles of $P'$, that are 
 the angles between the tangent spaces  $T_{p'}(Q'_i)$ and  $T_{p'}(Q'_j)$ at the points $p'$ on the $(n-2)$-faces  $Q'_{ij}=Q'_i\cap Q_j'$,  satisfy
 $$\angle_{ij}(P')\leq  \angle_{ij}(P),$$

$\bullet $ this  angle inequality is strict at some point, i.e.  there exits $p_0'\in Q'_{ij}$ in   some $Q'_{ij}$, such that 
 $$\angle(T_{p_0'}(Q'_i),T_{p_0'}(Q'_j))< \angle_{ij}(P).$$
\vspace {1mm}

 {\large \sf \textbf {F}$_1$. Question.} {\sl Are all extremal convex polyhedra $P$ are mean convexly extremal?} 
 
It is not even  known if the
 regular $3$-simplex is mean convexly extremal, but\vspace {1mm}
 
 \hspace {24mm} { \sl  the mean convex extremality of the $n$-cube}\vspace {1mm}
 
  \hspace {-6mm} follows by developing the cube  $P$ into  a complete  (orbi-covering) manifold $\hat P$  homeomorphic to $\mathbb R^n$ by reflecting $P$  in the faces, approximating the  natural continuous Riemannin metric metric on $\hat P$ by a smooth one  with $Sc\geq\varepsilon >0$  (see [28]\footnote {[28] M.Gromov {\sl Dirac and Plateau Billiards in Domains with Corners}, Cent. Eur. J. Math.
12(8) pp 1109-1156, (2014).}) and appealing to  gap  extremality of  $\mathbb R^n$ stated in  section 5.

And the same  argument   yields  (see [28]) the following \vspace {1mm}

 [{\Large $\ast$}]   {\sl  Let  a Riemannin metric $g$ on the   $n$-cube $P$   satisfy:

\hspace {6mm} $\ast_0$  $Sc(g)\geq 0$. 
  
\hspace {6mm}  $\ast_1$   $mean.curv_g (Q_i)\geq 0$,
   
 \hspace {6mm}    $\ast_2$ $\angle_{ij}(P,g)\leq \frac {\pi}{2}$.
    
   \hspace {-6mm} Then, necessarily,   $Sc(g)= 0$,   $mean.curv_g (Q_i)= 0$ and $\angle_{ij}(P,g)= \frac {\pi}{2}$.}
   
   \vspace {1mm}
   
  Probably, these equalities imply   that $P$ is {\it  isometric to  a Euclidean rectangular solid} but the approximation/smoothing   is no good  for proving this kind of rigidity.

  \vspace {1mm}
  
The main merit of    [{\Large $\ast$}]  is that it  provides a test  for  $Sc\geq 0$   in all  Riemannin manifolds $X$:
\vspace {1mm}

{\sl $Sc(X)\geq 0$ if and only if 
 no   cubical domain  $P\subset X$  satisfies
 $$[mean.curv_g (Q_i)>0]\& [\angle_{ij}(P,g)\leq \frac {\pi}{2}].$$}

 This suggests a possibility of  defining   $Sc(X)\geq 0$ for some singular spaces, $X$. e.g.
 for {\it Alexandrov spaces} with sectional curvatures bounded from below. 
 \vspace {1mm}

 {\large \sf \textbf {F}$_2$. Conjecture } {\sl All known (and expected) properties of Riemannian manifolds with  $Sc\geq 0$,  which have  no  "spin" attached to their formulations,  generalise to  Alexandrov's spaces.}\vspace {1mm}

For instance, most probably, \vspace {1mm}

  if an $n$-dimensional Alexandrov space $X$ with curvatures bounded from below has   $Sc> 0$ at all regular points $x\in X$, (or if the volumes of all   infinitesimally small balls in $X$ are  bounded by the volumes of such Euclidean balls) then\vspace {1mm}

  {\large \sf every continuous map from $X$ to a space $Y$ with $CAT(0)$ universal covering (i.e. an   Alexandrov's space   with  non-positive  sectional curvatures) contracts to an $(n-1)$-dimensional subset in $Y$.}\vspace {1mm}

If true, this would imply that (suitably defined) {\it harmonic maps} $X\to Y$ must necessarily  have {\it $(n-1)$-dimensional} images, which suggests a  (non-local?)  Weitzenboeck-Bochner type formula in this context and a definition of $Sc>0$ via spectral properties of small (large?) balls (cubes?) in $X$.

\section{References Listed.}

\hspace {4.3mm} \cite {[1]}    S.Alexander, V. Kapovitch, A.Petrunin,
{\sl Alexandrov 
geometry,}

\url{http://www.math.psu.edu/petrunin/ }\vspace {1mm}

 [2] David Hilbert, {\sl The Foundations of Physics}, (1915). \vspace {1mm}
 
[3]  A. Lichnerowicz,
{\sl Spineurs harmoniques.}
C. R. Acad. Sci. Paris, S\'erie A, 257 (1963), 7-9. \vspace {1mm}
 
[4] N. Hitchin,
{\sl Harmonic Spinors},
 Adv. in Math.
14
(1974), 1-55. \vspace {1mm}

[5]  S. Stolz. {\sl Simply   connected  manifolds   of positive  scalar  curvature},
  Ann.  of  Math.  (2) 136
  (1992),  511-540.\vspace {1mm}

[6]   J. Lohkamp,
{\sl  Metrics of negative Ricci curvature}, Annals of Mathematics, 140 (1994), 655-683.\vspace {1mm}

[7]  R. Schoen and S.-T. Yau, {\sl On the proof of the positive mass conjecture in general relativity,} Commun. Math. Phys. 65,  (1979). 45-76.\vspace {1mm}

[8]   M. Min-Oo, {\sl Scalar curvature rigidity of asymptotically hyperbolic spin manifolds}, Math. Ann. 285, 527?539 (1989). \vspace {1mm}

  [9]  L. Guth, {\sl  Notes on Gromov's systolic estimate}, Geom Dedicata (2006) 123:113-129. \vspace {1mm}

[10]   J Kazdan, F. Warner,
  {\sf Existence and Conformal Deformation of Metrics With Prescribed Gaussian and
Scalar Curvatures,} 
Annals of Mathematics,  101,\# 2. (1975), pp. 317-331.\vspace {1mm}

[11]  Yau, S.T., and Schoen, R.  {\sl On the Structure of Manifolds with positive Scalar Curvature}. Manuscripta mathematica 28 (1979): 159-184.\vspace {1mm}

 [12] J. Lohkamp,
{\sl The Higher Dimensional Positive Mass Theorem II}, (2016)
 	arXiv:1612.07505.\vspace {1mm}

[13] R. Schoen, S.T. Yau, {\sl Positive Scalar Curvature and Minimal Hypersurface Singularities},
 (2017) 	arXiv:1704.05490. \vspace {1mm}

[14] M. Gromov, H. B. Lawson, Jr., {\sl Positive scalar curvature and the Dirac operator on complete Riemannian manifolds},  Publ. Math. IH\'ES 58 (1983), 295-408. \vspace {1mm}

[15] M. Gromov, {\sl Positive curvature, macroscopic dimension, spectral gaps and higher signatures,}  in Proc of 1993 Conf. in Honor of of the Eightieth Birthday of I. M. Gelfand, 
 Functional Analysis on the Eve of the 21st Century: Volume I  Progress in Mathematics,  (1996) pp 1-213, Vol. 132, \vspace {1mm}

[16] S. Markvorsen, M. Min-Oo,  {\sl Global Riemannian Geometry: Curvature and Topology}, 2012
Birkh\"auser.\vspace {1mm}

 [17] L. Guth, {\sl Metaphors in systolic geometry}. In: Proceedings of the International Congress of Mathematicians.  2010,  Volume II, pp. 745-768.\vspace {1mm}

[18] L. Guth,{\sl  Volumes of balls in Riemannian manifolds and Uryson width}.
Journal of Topology and Analysis,
Vol. 09, No. 02, pp. 195-219 (2017).\vspace {1mm}

[19] D.Bolotov, A. Dranishnikov.  {\sl On Gromov's conjecture for totally non-spin manifolds}, (2015)
arXiv:1402.4510v6.\vspace {1mm}

[20] M. Marcinkowski, {\sl Gromov positive scalar curvature conjecture and rationally inessential macroscopically large manifolds}.
 Journal of topology  9, 1; 105-116  (2016).
 Oxford University Press. \vspace {1mm}

[21] J. Rosenberg,  {\sl Manifolds of positive scalar curvature: a progress report}, in: Surveys on Differential Geometry,
vol. XI: Metric and Comparison Geometry, International Press 2007.\vspace {1mm}

[22] M. Gromov, {\sl Morse Spectra, Homology Measures, Spaces of Cycles and Parametric Packing Problems.}

\url {www.ihes.fr/~gromov/PDF/Morse-Spectra-April16-2015-.pdf} \vspace {1mm}

 [23]  M. Llarull,   {{\sl Sharp estimates and the Dirac Operator},

\hspace {2mm} Math. Ann.
310 (1998), 55-71. \vspace {1mm}

 [24]  M. Min-Oo, {\sl Scalar Curvature Rigidity of Certain Symmetric Spaces}, 
 Geometry, Topology and Dynamics
(Montreal, PQ, 1995), CRM Proc. 
  Lecture Notes, 15, Amer. Math. Soc., Providence, RI, 1998, pp. 127-136.\vspace {1mm}

 [25]  S. Goette and U. Semmelmann, {\sl Scalar curvature estimates for 
 
  compact symmetric spaces}, 
 Differential Geom. Appl. 16(1), (2002) 65-78. \vspace {1mm}

[26]  M.Listing, {\sl   The  Scalar curvature on compact symmetric spaces},
 arXiv:1007.1832, 2010 - arxiv.org. \vspace {1mm}

[27] M.Gromov,
{\sl Metric Inequalities with Scalar Curvature.}

\url
{http://www.ihes.fr/~gromov/PDF/Inequalities-July\%202017.pdf } \vspace {1mm}

 [28] {M.Gromov {\sl Dirac and Plateau Billiards in Domains with Corners}, Cent. Eur. J. Math.
12(8) pp 1109-1156, (2014).

\end{document}